\documentclass[9pt]{revtex4}
 \usepackage{verbatim}
\usepackage{graphicx}% Include figure file
\usepackage{float}
\usepackage{amssymb}
\usepackage[utf8]{inputenc}
\usepackage[english,french]{babel}
\usepackage[intlimits]{amsmath}
\usepackage{amsfonts}
\usepackage{bm}%
\usepackage{hyperref}

 \newtheorem{thm}{Theorem}
 
 \newtheorem{lem}[thm]{Lemma}
 \newtheorem{prop}[thm]{Proposition}
\usepackage{amsmath,amscd}
\usepackage{amsfonts}
\usepackage{bm}
\usepackage{hyperref}

 \newcommand{\R}{\mathbb{R}}

 \newcommand{\tr}{\mathcal{T}}
 \newcommand{\tp}{\widehat{\mathcal{T}}}

 \newcommand{\Pro}{\mathbb{P}^\mathrm{st}}

 \newcommand{\abs}[1]{\left\vert#1\right\vert}

\begin{document}
\title{TASEP and Planar Binary Trees: A Combinatorial Approach}
\author{Xiangyu Cao}
\affiliation{Institut de Physique Théorique, CEA Saclay / Campus Paris--Saclay, France.}

\begin{abstract}
We give a simple explanation why the stationary state of the 1D TASEP model with open boundaries is related to the Catalan numbers. Our construction is based on planar binary trees and provides a combinatorial solution of the stationary state. We discuss the relation to permutation/Catalan tableaux approaches in an appendix.
\end{abstract}
\maketitle

\section{Introduction} 
The Totally Asymmetric Simple Exclusion Process with open boundary condition(open--boundary TASEP) is one of the first interacting driven lattice gas models whose \textit{stationary state} has been obtained exactly \cite{zeilberger,derrida,derrida2}; The stationary state of a driven, out--of--equilibrium particle system is unknown, contrary to a system at equilibrium for which one can write down the Boltzmann distribution. In fact, the open--boundary TASEP belongs to a large family of stochastic models whose stationary state can be calculated by algebraic--integrability techniques, notably the matrix product Ansatz \cite{matrix}. On the other hand, the stationary state of these models has rich combinatorial structures. In particular the stationary weights of the open--boundary TASEP are known to be closely related to Catalan numbers. This relation has been interpreted using Catalan tableaux(see Viennot et. al. \cite{viennot}) or more generally using permutation tableaux(\cite{corteel2,corteel3}). Duchi and Schaeffer(\cite{duchi}) constructed a TASEP model with one hidden line to derive the stationary measure combinatorially. Corteel et. al. managed to do so using lattice paths(\cite{corteel}) or permutation tableaux (\cite{corteel4}) for more general partially asymmetric exclusion processes. The relation to tableaux combinatorics can also lead to determinant formulae for probabilities (\cite{mandel}).

Here we would like to record another simple combinatorial approach to the open--boundary TASEP stationary state, based on binary trees. Although binary trees are also used in Angel's combinatorial approach (\cite{angel}) to the TASEP with two species of particles(for a definition see \cite{matrix}, Sect. 4), our approach is more closely related to permutation/Catalan tableaux approaches, and particularly close to that of \cite{corteel4}. The appendix \ref{app:tableaux} will be devoted to discussing this relation while the rest of this note will be self--contained.
\vspace{0.2cm}

\paragraph*{\bf{Setup and main result.}}
In what follows let $n>0$ be a positive integer. The $n$-site open--boundary TASEP is a continuous time Markov chain whose configuration space is $\{\bullet,\circ\}^n$. A state $\mathcal{C}$ therein is a string of length $n$ of $\bullet / \circ$'s and describes the particle occupancy of $n$ sites. Its dynamics is determined by the following \textit{transition rates}
$$ W(\mathcal{C} \rightarrow \mathcal{C}') = \begin{cases} 
\alpha, & \text{ if } \mathcal{C} = \circ\mathcal{A}, \mathcal{C}' = \bullet\mathcal{A}; \\
\beta,  & \text{ if } \mathcal{C} = \mathcal{A}\bullet, \mathcal{C}' = \mathcal{A} \circ; \\
1,      & \text{ if } \mathcal{C} = \mathcal{A} \bullet\circ\mathcal{A}'; \mathcal{C}' = \mathcal{A} \circ\bullet \mathcal{A}'; \\
0,      & \text{ otherwise. }
\end{cases} $$
where $\alpha, \beta\in (0,1]$ are parameters and $\mathcal{A},\mathcal{A}'$ represent arbitrary $\circ$/$\bullet$ strings. A nontrivial question is to determine the \textit{stationary measure} of such a system for any $n$. This amounts to finding $\Pro_n\{\mathcal{C}\}$ such that for any configuration $\mathcal{C}$,
$$
\Pro_n\{\mathcal{C}\} \sum_{\mathcal{C}'} W(\mathcal{C} \rightarrow\mathcal{C}') = 
 \sum_{\mathcal{C}'} \Pro_n\{\mathcal{C}'\} W(\mathcal{C}' \rightarrow{C}).
$$

\begin{figure}[h]
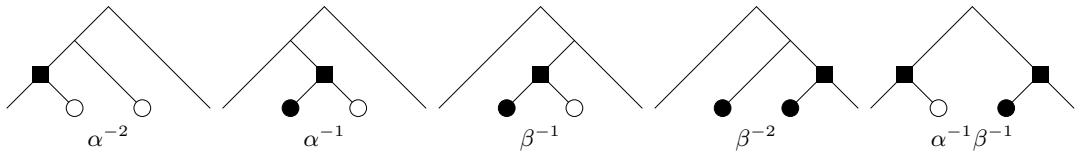

\begin{tabular}{ccccc}
\includegraphics[scale=0.45]{trees.1} &
\includegraphics[scale=0.45]{trees.2} &
\includegraphics[scale=0.45]{trees.4} &
\includegraphics[scale=0.45]{trees.3} &
\includegraphics[scale=0.45]{trees.5}  \\ 
$\alpha^{-2}$ & $\alpha^{-1}$ & $\beta^{-1}$ & $\beta^{-2}$ & $\alpha^{-1}\beta^{-1}$ 
\end{tabular}
\caption{Elements of $\mathcal{T}_2$, their weights and their induced configuration. The marks $\blacksquare$ will be explained in the text.}
\label{trees2}
\end{figure}

Let us denote $\tr_n$ the set of plane binary trees  with $(n+2)$-endpoints. Their cardinality is given by the Catalan numbers  $\abs{\tr_n} = C_{n+1} = \frac{1}{n+2} \binom{2n+2}{n+1}$ (\cite{stanley}). For example the $5$ elements of $\mathcal{T}_2$ shown in figure \ref{trees2}. 

On $\tr_n$ we shall define a \textit{reduction map} $R: \tr_n \rightarrow \{\bullet,\circ\}^n$ and a \textit{weight function} $\mu: \tr_n \rightarrow \R^+$ as follows:
\begin{itemize}
\item[*] For $T \in \tr_n$, $R(T) := (t_1, \dots, t_n)$ (called its \textit{reduced configuration}) and $t_k$ is determined by the $(k+1)$-th endpoint of $T$ counting from left: $t_k = \bullet$ if it is a left descendent(its edge is of up-right direction) and $t_k = \circ$ if it is a right descendent(its edge is of up-left direction). The leftmost and rightmost endpoints are not concerned. For $n=2$, from left to right in figure \ref{trees2}, $R(T)$'s are $\circ\circ,\bullet\circ,\bullet\circ, \bullet\bullet$ and $\circ\bullet$. $R$ is surjective but in general not one--to--one: For example $R^{-1}\{\bullet\circ\}$ has $2$ elements.  
\item[*] $\mu(T) := \alpha^{-l(T)} \beta^{-r(T)}$ where $l(T)$(respectively $r(T)$) is the number of vertices on the path between the leftmost(rightmost, respectively) endpoint and the root. For example, $\mu\begin{pmatrix} \includegraphics[scale=0.15]{trees3.5}\end{pmatrix} = \alpha^{-2}$ and weights of elements in $\tr_2$ are shown in \ref{trees2}.
\end{itemize}

\begin{prop}\label{prop:main2}
The stationary measure of the $n$-TASEP with parameter $\alpha,\beta$ is given by 
$$ \Pro_n\{\mathcal{C}\} = (Z_n)^{-1} \sum_{T\in R^{-1} \{\mathcal{C}\}} \mu(T) $$
where $\mu(T) = \alpha^{-l(T)} \beta^{-r(T)} $ where $l(T)$($r(T)$) is the number of vertices between the leftmost(rightmost) endpoint of $T$ and the its root. $Z_n = \sum_{T\in\tr_n} \mu(T)$ is the normalisation factor.
\end{prop}

The rest of the note is devoted to a combinatorial proof of the above proposition. For the sake of transparency, we shall prove it first for the special case of $\alpha = \beta = 1$ in Sect. \ref{sec:special}. The general case is considered in Sect. \ref{sec:general}.

\section{Case of $\alpha = \beta = 1$}\label{sec:special}
In this section we focus on the case of $\alpha = \beta = 1$, \textit{i.e.}, and consider the following restriction of Proposition \ref{prop:main2}:
\begin{prop}\label{prop:main}
The stationary probability of the $n$-TASEP with $\alpha = \beta = 1$ is induced by the uniform measure on $\tr_n$
$$
\Pro_n\{\mathcal{C}\} = \frac{\abs{ R^{-1} \{\mathcal{C}\} }}{\abs{\tr_n}}. 
$$
(Here and below, $\abs{X}$ denotes the cardinal of the set $X$.)
In other words, for any state $\mathcal{C}\in \{\bullet, \circ\}^n$, we have the flux balance equation
\begin{equation}\label{eqn:bilan}
 \abs{ R^{-1} \{\mathcal{C}\}} \sum_{\mathcal{C}'} W(\mathcal{C} \rightarrow\mathcal{C}') = 
 \sum_{\mathcal{C}'} \abs{ R^{-1} \{\mathcal{C}'\}} W(\mathcal{C}' \rightarrow\mathcal{C}).
\end{equation}
\end{prop}
Since $\alpha=\beta=1$, $W(\mathcal{C} \rightarrow\mathcal{C}')$ is either $0$ or $1$, so both sides of equation \ref{eqn:bilan} are integers, which suggests a combinatorial proof. Indeed, we shall interpret either side of (\ref{eqn:bilan}) as the size of a set, and construct a bijection between them.

On a plane binary tree, we define a \textit{last branching vertex} to be a non--endpoint vertex whose descendents are both endpoints. The last branching vertices were all marked with $\blacksquare$'s in the figure \ref{trees2}. We define a \textit{marked} plane binary tree is a plane binary tree with \textit{only one} of its last branching vertices marked with a $\blacksquare$. We denote $\tp_n$ the set of marked binary trees of $(n + 2)$ end-points and let $f:\tp_n \rightarrow \tr_n$ be the map that forgets the marking. For $n > 1$, we have $\abs{\tp_n} > \abs{\tr_n}$ since $f$ is subjective but not one--to--one. For example, $\tp_2$ has $6$ elements as listed below:
\begin{figure}[H]
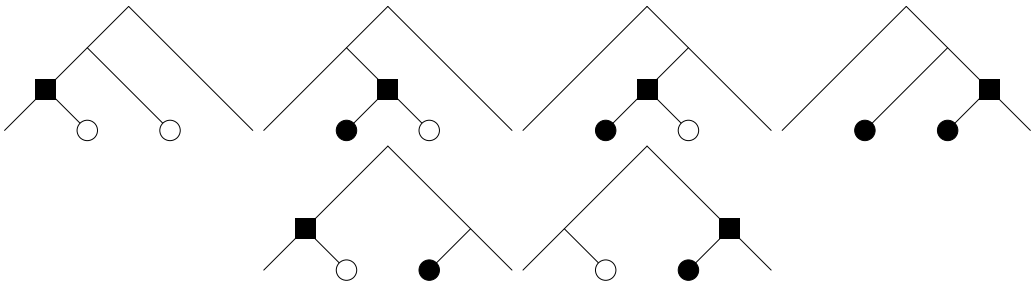

\center
\includegraphics[scale=0.55]{trees.1}
\includegraphics[scale=0.55]{trees.2}
\includegraphics[scale=0.55]{trees.4}
\includegraphics[scale=0.55]{trees.3}\\
\includegraphics[scale=0.55]{trees.6}
\includegraphics[scale=0.55]{trees.7}
\caption{Elements of $\tp_2$. Note that the last two elements in figure \ref{mtrees2} are mapped by $f$ to the same element of $\tr_2$(namely, the last one in figure \ref{trees2}).}\label{mtrees2}
\end{figure}
$\tp_3$ has $20$ elements(while $\tr_3$ has $14$), we draw them in the appendix \ref{app:t3}.(In general, one can show that $\tr_n$ has $\binom{2n}{n}$ elements. But we will not need this result for the following)

Now we make a crucial observation that follows from the above definitions: The last branching vertices of a tree $T$ correspond to the active bonds(places where a transition to another configuration is possible) of its reduced configuration $R(T)$. As a consequence, we have the following
\begin{lem}
Suppose $\alpha = \beta = 1$. For any $\mathcal{C} \in \{\bullet,\circ\}^n$ and $T\in R^{-1}\{\mathcal{C}\}$, we have $\sum_{\mathcal{C}'} W(\mathcal{C} \rightarrow\mathcal{C}') = \abs{f^{-1}\{T\}}$. As a consequence, for any $\mathcal{C} \in \{\bullet,\circ\}^n$, 
\begin{equation}\label{eq:lhs}
\abs{ R^{-1} \{\mathcal{C}\}} \sum_{\mathcal{C}'} W(\mathcal{C} \rightarrow\mathcal{C}') = \abs{ (Rf)^{-1}\{\mathcal{C}\}}.
\end{equation}
\end{lem}

Thus we managed to see the left hand side of \ref{eqn:bilan} as the size of a set. Next, we shall describe a map $ \pi: \tp_n \rightarrow \tp_n $ which we will prove to be a bijection, and such that the right hand side of \ref{eqn:bilan} is the size of the set $ (Rf\pi)^{-1}\{\mathcal{C}\}$, which is bijected to 
$ f^{-1}R^{-1}\{\mathcal{C}\}$ by $\pi$. This would conclude the proof of Proposition \ref{prop:main}.

Consider a marked tree $\widehat{T}$. Since the tree root is never a last branching vertex (for $n>0$), its marked last branching vertex $v$ is either a left descendent or a right descendent (of its ascendent). Suppose for the moment that $v$ is a right descendent. To obtain $\pi(\widehat{T})$, we remove the segment linking $v$ and its left descendent, glue it to the next right-descendent endpoint to its right and mark the new vertex. For example:
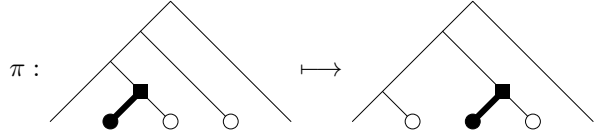
\begin{center}
$\pi:$ $\begin{matrix}\includegraphics[scale=0.4]{trees3.23}\end{matrix}$ $\longmapsto$ $\begin{matrix}\includegraphics[scale=0.4]{trees3.24} \end{matrix}$
\end{center}
To illustrate, the moved segment is drawn in bold. 

In the special case where $v$ has the rightmost endpoint of the tree as descendent, the above operation cannot be done. In this case, we glue the removed segment to the rightmost left-descendent endpoint of the tree and mark the new vertex, for example 
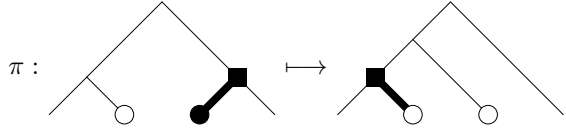
\begin{center}
$\pi:$ $\begin{matrix}\includegraphics[scale=0.5]{trees.8}\end{matrix}$ $\longmapsto$ $\begin{matrix}\includegraphics[scale=0.5]{trees.9}\end{matrix}$
\end{center}
Note that in this second (boundary) case the moved segment changes direction as opposed to the first case. 

So we have defined $\pi$ when the marked last branching vertex is a right descendent. For the alternative possibility, we just need to swap the words ``left'' and ``right'' in the recipe described above (Equivalently, one may say that we complete the definition of $\pi$ in the unique way so that it preserves the left--right symmetry). We illustrate below $\pi$ for $n = 2$ (see appendix \ref{app:t3} for $n = 3$).
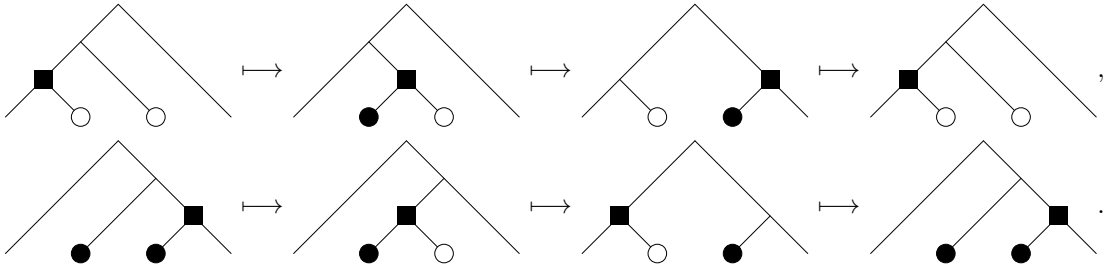
\begin{center}
$\begin{matrix}\includegraphics[scale=0.5]{trees.1}\end{matrix}$ $\longmapsto$ $\begin{matrix}\includegraphics[scale=0.5]{trees.2}\end{matrix}$ $\longmapsto$ $\begin{matrix}\includegraphics[scale=0.5]{trees.7}\end{matrix}$ $\longmapsto$ $\begin{matrix}\includegraphics[scale=0.5]{trees.1}\end{matrix}$, \\
$\begin{matrix}\includegraphics[scale=0.5]{trees.3}\end{matrix}$ $\longmapsto$ $\begin{matrix}\includegraphics[scale=0.5]{trees.4}\end{matrix}$ $\longmapsto$ $\begin{matrix}\includegraphics[scale=0.5]{trees.6}\end{matrix}$ $\longmapsto$ $\begin{matrix}\includegraphics[scale=0.5]{trees.3}\end{matrix}$. 
\end{center}

One sees that they are permutations(bijections). This is true in general:
\begin{lem}\label{lem:bij}
$\pi: \tp_n \rightarrow \tp_n$ is a bijection.  
\end{lem}
In fact we can describe its inverse $\sigma$ in a similar manner. To obtain $\sigma(\widehat{T})$ from $\widehat{T}$, assuming that the marked last branching vertex $v$ is a right(left) descendent, remove the segment linking $v$ and its left(right) descendent, glue it to the next right(left)-descendent endpoint to its left(right); when it is impossible to do so, move it to the leftmost(rightmost) endpoint (changing its direction). In either case, mark the new vertex. This construction guarantees that $\sigma$ and $\pi$ are inverse to each other. 

Now let us take a pair $\widehat{T} = \pi (\widehat{T}')$ and consider their respective reduced states $ \mathcal{C} = Rf(\widehat{T})$ and $\mathcal{C}'= Rf(\widehat{T}')$. One can observe by construction of $\pi$ that $W(\mathcal{C}'\rightarrow \mathcal{C}) \neq 0$(in fact equal to $1$ provided $\alpha = \beta = 1$) for any such pair $\widehat{T}$ and $\widehat{T}'$. 

Conversely, for any pair of states $\mathcal{C},\mathcal{C}'$ such that $W(\mathcal{C}' \rightarrow \mathcal{C}) \neq 0$ ($=1$ in the case $\alpha = \beta = 1$), there are exactly $\abs{R_n^{-1}\{\mathcal{C}'\}}$ marked trees $\widehat{T}'$ such that $Rf\pi(\widehat{T}') = \mathcal{C}$. In fact, for each $T' \in R^{-1} \{\mathcal{C}'\}$, there is one and only one $\widehat{T}'$ such that $f(\widehat{T}') = T'$ and that $Rf\pi(\widehat{T}') = \mathcal{C}$: It is the one obtained by marking the last branching vertex corresponding to the active bond involved in the transition $\mathcal{C}' \rightarrow \mathcal{C}$. 

We have thus shown the following:
\begin{lem}
For any $\mathcal{C},\mathcal{C}' \in \{0,1\}^{n}$, the set $\mathcal{F}(\mathcal{C}' \rightarrow \mathcal{C})$ defined as
\begin{equation}
\mathcal{F}(\mathcal{C}' \rightarrow \mathcal{C}) := (Rf)^{-1}\{\mathcal{C}'\} \bigcap (Rf\pi)^{-1}\{\mathcal{C}\} 
\end{equation}
is non-empty if and only if $W(\mathcal{C}' \rightarrow \mathcal{C}) \neq 0$. In this case its cardinal is equal to that of $R^{-1}\{\mathcal{C}'\}$. Thus, when $\alpha = \beta = 1$ we have
\begin{equation} \label{eq:flux}
\abs{\mathcal{F}(\mathcal{C}' \rightarrow \mathcal{C})}  =  \abs{R^{-1}\{\mathcal{C}'\}} W(\mathcal{C}' \rightarrow \mathcal{C}).
\end{equation}
\end{lem}
Summing the above equation over $\mathcal{C}'$, we have for any $\mathcal{C} \in \{ \bullet, \circ\}^{n}$
\begin{equation}\label{eqn:rhs}
\abs{(Rf\pi)^{-1}\{\mathcal{C}\}} = \sum_{\mathcal{C}'} \abs{ (Rf)^{-1}\{\mathcal{C}'\} \bigcap (Rf\pi)^{-1}\{\mathcal{C}\} } = 
\sum_{\mathcal{C}'} \abs{ R^{-1} \{\mathcal{C}'\}} W(\mathcal{C}' \rightarrow\mathcal{C})
\end{equation}
which is promised before. Now combining equations (\ref{eq:lhs}), the bijectiveness of $\pi$ (lemma \ref{lem:bij}) and the equation \ref{eqn:rhs}) gives 
$$ \abs{ R^{-1} \{\mathcal{C}\}} \sum_{\mathcal{C}'} W(\mathcal{C} \rightarrow\mathcal{C}') = \abs{(Rf)^{-1}\{\mathcal{C}\}} = 
\abs{(Rf\pi)^{-1}\{\mathcal{C}\}} =  \sum_{\mathcal{C}'} \abs{ R^{-1} \{\mathcal{C}'\}} W(\mathcal{C}' \rightarrow\mathcal{C}),$$
which proves Proposition \ref{prop:main}. Notice that summing eq.(\ref{eq:flux}) over $\mathcal{C}$ gives eq. (\ref{eq:lhs}), so it is in fact the only ingredient for proving Proposition \ref{prop:main} besides the bijectiveness of $\pi$.

\section{General $\alpha,\beta$}\label{sec:general}
Now that combinatorial constructions are independent of the value of $\alpha$ and $\beta$, to adapt our proof to arbitrary values of $\alpha$ and $\beta$, we need to assign weights to marked trees so that a weighted version of (\ref{eq:flux}) holds. 
We define $\hat{\mu}:\tp_n \rightarrow \R^+$ as $\hat{\mu}(\widehat{T}) = \alpha^{-l'(\widehat{T})} \beta^{-r'(\widehat{T})}$ where $l'(\widehat{T})$(respectively $r'(\widehat{T})$) counts the same vertices as do $l(T)$($r(T)$ respectively) \textit{except} the marked one when it is there: For example when $\widehat{T} = \begin{matrix} \includegraphics[scale=0.3]{trees3.20} \end{matrix}$, $\hat{\mu} \left(\widehat{T}  \right) = \alpha^{-1}\beta^{-1}$ while $\mu(f(\widehat{T})) = \alpha^{-1}\beta^{-2}$. $\mu$ ($\hat{\mu}$ respectively) defines a measure on (the power set of) $\tr_n$ (respectively $\tp_n$) in the obvious way and we shall denote the measure again by $\mu$($\hat{\mu}$ respectively). 

To show Proposition \ref{prop:main2}, all we need is the following two properties of $\mu$ and $\hat{\mu}$, which can be readily checked :
\begin{enumerate}
\item[(a)] For any $\widehat{T} \in \mathcal{F}(\mathcal{C}' \rightarrow \mathcal{C})$, 
$\hat{\mu}(\widehat{T}) = \mu(f(\widehat{T})) W(\mathcal{C}' \rightarrow \mathcal{C})$.
\item[(b)] For any $\widehat{T} \in \tp_n$, $\hat{\mu}(\widehat{T}) = \hat{\mu}(\pi(\widehat{T}))$.
\end{enumerate}
Property (a) implies the following variant of (\ref{eq:flux})
%$$
%\sum_{\widehat{T}\in\mathcal{F}(\mathcal{C}' \rightarrow \mathcal{C})} \hat{\mu}(\widehat{T}) = \sum_{T\in R^{-1}\{\mathcal{C}'\}} \mu(T)W(\mathcal{C}'\rightarrow \mathcal{C}).
%$$
\begin{equation}\label{eqn:flux2}
\hat{\mu} (\mathcal{F}(\mathcal{C}' \rightarrow \mathcal{C}))  = \mu(R^{-1} \{\mathcal{C}'\})W(\mathcal{C}'\rightarrow \mathcal{C}).
\end{equation}
 Summing it over $\mathcal{C}'$ gives 
$$
 \hat{\mu}((Rf\pi)^{-1}\{\mathcal{C}\}) = \sum_{\mathcal{C}'}  \mu(R^{-1} \{\mathcal{C}'\}) W(\mathcal{C}' \rightarrow\mathcal{C})
$$ which generalises eq. (\ref{eqn:rhs}), while summing (\ref{eqn:flux2}) over $\mathcal{C}$(and switching variable names $\mathcal{C}$ and $\mathcal{C}'$) gives 
$$ \hat{\mu}((Rf)^{-1}\{\mathcal{C}\}) = \sum_{\mathcal{C}'} \mu(R^{-1} \{\mathcal{C}\}) W(\mathcal{C} \rightarrow\mathcal{C}')$$ which generalises eq. (\ref{eq:lhs}). Now property (b) implies 
$$ \hat{\mu}((Rf)^{-1}\{\mathcal{C}\}) = \hat{\mu}((Rf\pi)^{-1}\{\mathcal{C}\} ),$$ hence 
$$ \sum_{\mathcal{C}'}  \mu(R^{-1} \{\mathcal{C}\} ) W(\mathcal{C} \rightarrow\mathcal{C}') =  \sum_{\mathcal{C}'} \mu(R^{-1} \{\mathcal{C}'\}) W(\mathcal{C}' \rightarrow\mathcal{C}). $$
That is, the stationary measure of the $n$-site TASEP is given by 
$ \Pro_n\{\mathcal{C}\} = (Z_n)^{-1} \mu(R^{-1} \{\mathcal{C}\})$, where $Z_n = \mu(\tr_n)$,  as claimed in Proposition \ref{prop:main2}.

\section{Conclusion}
In this note we have considered a combinatorial solution of the open boundary TASEP stationary state in terms of planar binary trees. We can summarise our approach as two points. First, the (nontrivial) stationary measure on the configuration space is induced by a simple measure on the space of plane binary trees (in the case of $\alpha = \beta = 1$ the latter measure is uniform). Second, the balance equation is given a bijective proof(at least in the $\alpha = \beta = 1$ case): every term therein is interpreted as the cardinal of a set, and equality comes from the construction of a bijection. 

Let us mention briefly some models whose stationary state can be obtained by a variant of the approach presented above, leaving details for future work. One can adapt this approach to the periodic boundary condition TASEP with one second class particle(for an introduction see \cite{matrix}, Sect. 4) with only a slight modification of $\pi$. Essentially the same construction works when there are $k$ second class particles: For this we need to consider $k$-tuples of planar binary trees. The latter construction works also for a special case of open boundary two--species TASEP where the second class particles are neither injected nor extracted, as first considered in \cite{2species}. An interesting question is whether it is possible to generalise our approach to the multi--species TASEP, whose stationary state has been obtained combinatorially by Martin and Ferrari in \cite{multi}, but with a queueing--theory--inspired construction.

\appendix
\section{Elements of $\tp_3$ and the permutation $\pi$}\label{app:t3}
We first enumerate all elements of $\tp_3$. They are organised in the following manner
\begin{itemize}
\item A configuration $\mathcal{C}\in \{\bullet, \circ\}^3$ is followed by the elements of $f^{-1}R^{-1}\{\mathcal{C}\}$.
\item Elements having the same image under $f$ are put in a parenthesis. 
\end{itemize}
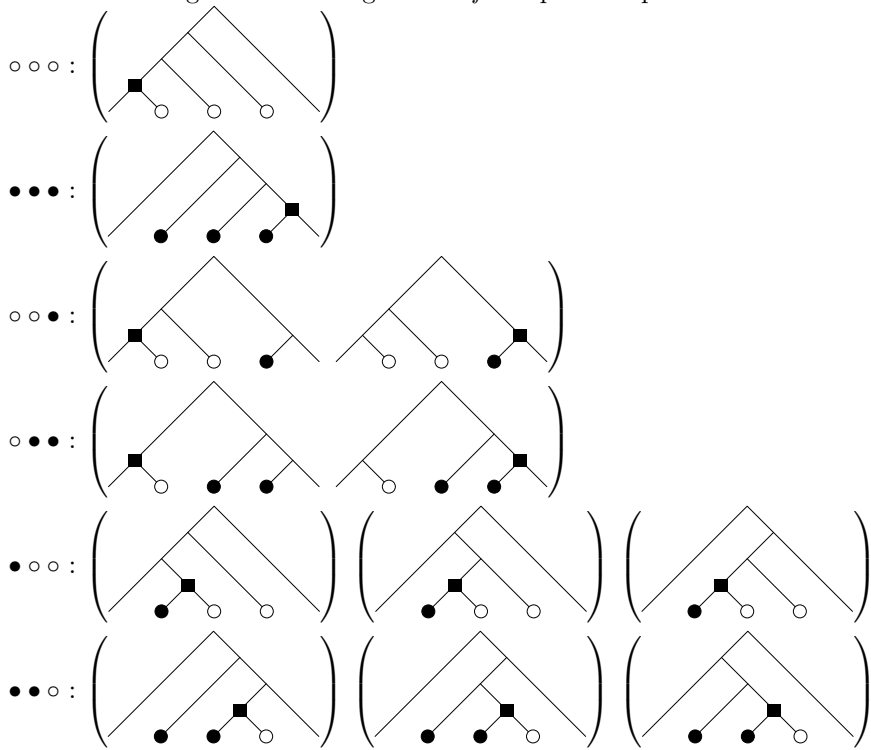
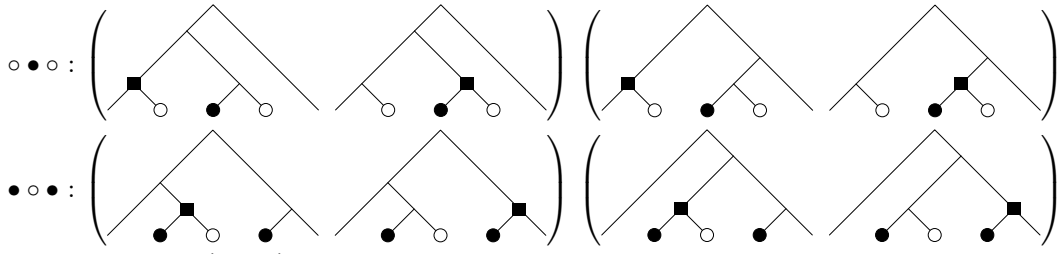
\begin{itemize}
\item[] $\circ\circ\circ$:  $\begin{pmatrix} \includegraphics[scale=0.35]{trees3.1}\end{pmatrix}$
\item[] $\bullet\bullet\bullet$:  $\begin{pmatrix} \includegraphics[scale=0.35]{trees3.8} \end{pmatrix}$
\item[] $\circ\circ\bullet$:    $\begin{pmatrix} \includegraphics[scale=0.35]{trees3.19}    & \includegraphics[scale=0.35]{trees3.21} \end{pmatrix}$ 
\item[] $\circ\bullet\bullet$:   $\begin{pmatrix} \includegraphics[scale=0.35]{trees3.22}  & \includegraphics[scale=0.35]{trees3.20}   \end{pmatrix}$                       
\item[] $\bullet\circ\circ$:  $\begin{pmatrix}  \includegraphics[scale=0.35]{trees3.2}\end{pmatrix}$  $\begin{pmatrix} \includegraphics[scale=0.35]{trees3.4}\end{pmatrix}$  $\begin{pmatrix} \includegraphics[scale=0.35]{trees3.10}  \end{pmatrix}$
\item[] $\bullet\bullet\circ$:  $\begin{pmatrix} \includegraphics[scale=0.35]{trees3.9}\end{pmatrix}$  $\begin{pmatrix}  \includegraphics[scale=0.35]{trees3.11}\end{pmatrix}$ $\begin{pmatrix} \includegraphics[scale=0.35]{trees3.3}\end{pmatrix}$
\item[] $\circ\bullet\circ$:  $\begin{pmatrix}\includegraphics[scale=0.35]{trees3.6}&
                    \includegraphics[scale=0.35]{trees3.7}\end{pmatrix}$
                     $\begin{pmatrix}\includegraphics[scale=0.35]{trees3.18} &
                     \includegraphics[scale=0.35]{trees3.16} \end{pmatrix}$ 
\item[] $\bullet\circ\bullet$:  $\begin{pmatrix} \includegraphics[scale=0.35]{trees3.15}
                   & \includegraphics[scale=0.35]{trees3.17} \end{pmatrix}$
                    $\begin{pmatrix} \includegraphics[scale=0.35]{trees3.14} &
                     \includegraphics[scale=0.35]{trees3.13} \end{pmatrix}$
\end{itemize}
Next we describe $\pi: \tp_3 \rightarrow \tp_3$. It is a permutation, and we depict each of its cycles in one row: $\pi$ maps each element to the one to its right, and the rightmost element to the leftmost.
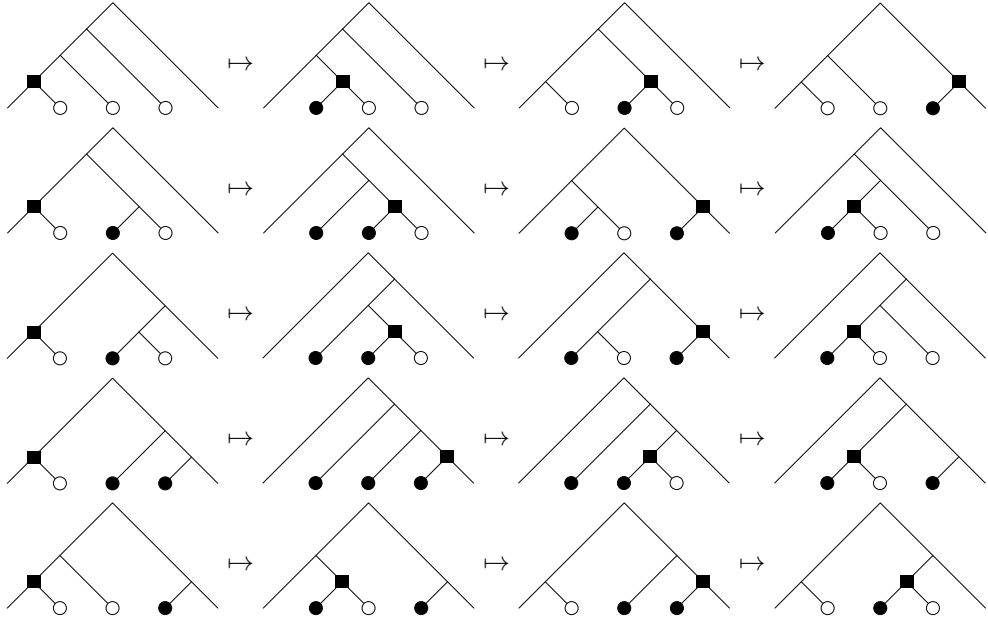
\begin{center}
$\begin{matrix}\includegraphics[scale=0.35]{trees3.1}\end{matrix}$ $\mapsto$ $\begin{matrix}\includegraphics[scale=0.35]{trees3.2}\end{matrix}$ $\mapsto$ $\begin{matrix}\includegraphics[scale=0.35]{trees3.7} \end{matrix}$ $\mapsto$ $\begin{matrix}\includegraphics[scale=0.35]{trees3.21} \end{matrix}$ \\ 
$\begin{matrix}\includegraphics[scale=0.35]{trees3.6}\end{matrix}$ $\mapsto$ $\begin{matrix}\includegraphics[scale=0.35]{trees3.3}\end{matrix}$ $\mapsto$ $\begin{matrix}\includegraphics[scale=0.35]{trees3.17} \end{matrix}$ $\mapsto$ $\begin{matrix}\includegraphics[scale=0.35]{trees3.4}\end{matrix}$ \\
$\begin{matrix}\includegraphics[scale=0.35]{trees3.18}\end{matrix}$ $\mapsto$ $\begin{matrix}\includegraphics[scale=0.35]{trees3.11}\end{matrix}$ $\mapsto$ $\begin{matrix}\includegraphics[scale=0.35]{trees3.13}\end{matrix}$ $\mapsto$ $\begin{matrix}\includegraphics[scale=0.35]{trees3.10}\end{matrix}$ \\
$\begin{matrix} \includegraphics[scale=0.35]{trees3.22}\end{matrix}$ $\mapsto$ $\begin{matrix}\includegraphics[scale=0.35]{trees3.8}\end{matrix}$ $\mapsto$ $\begin{matrix}\includegraphics[scale=0.35]{trees3.9}\end{matrix}$ $\mapsto$ $\begin{matrix}\includegraphics[scale=0.35]{trees3.14}\end{matrix}$\\
$\begin{matrix}\includegraphics[scale=0.35]{trees3.19}\end{matrix}$ $\mapsto$ $\begin{matrix}\includegraphics[scale=0.35]{trees3.15}\end{matrix}$ $\mapsto$ $\begin{matrix}\includegraphics[scale=0.35]{trees3.20}\end{matrix}$ $\mapsto$ $\begin{matrix}\includegraphics[scale=0.35]{trees3.16} \end{matrix}$\\ 
\end{center}
%Remark that, in general, the cycles of $\pi: \tp_n \rightarrow \tp_n$ are all of length $n+1$ and the set of cycles is bijectively mapped onto $ \abs{\tr_{n-1}}$ by removing the ``moving'' endpoint: This is closely related to the calculation average current of the open boundary TASEP.
\section{Relation with Catalan Tableaux}\label{app:tableaux}
This appendix comments briefly the relation of our planar binary tree construction with the combinatorial solution of open--boundary TASEP(and more general partially asymmetric processes) stationary measure by Corteel and Williams (\cite{corteel4}). We shall construct a simple correspondence between planar binary trees and Catalan tableaux, on which the above cited work is based(when restricted to the TASEP). 

Be definition, \textit{Catalan tableaux} are Young diagrams filled with $0$ and $1$'s($0/1$ tableaux) fulfilling two properties: 
\begin{enumerate}
\item On each column there is exactly one $1$, with the exception of the leftmost column where there is no $1$.
\item There are no $0$ that has a $1$ above it in the same column and a $1$ left to it in the same row. (The Young diagram is drawn in the down--right quadrant)
\end{enumerate} 
The \textit{index} of of a Catalan tableaux is defined as $\lambda_1 + \lambda_1^t - 1$, where $\lambda = (\lambda_1 \geq \lambda_2 \geq \dots)$ is the partition of the underlying Young diagram and $\lambda^t$ the transpose. For example, $\begin{matrix} \includegraphics[scale=0.6]{trees3.30} \end{matrix}$ is a Catalan tableaux of index $5$ while $\begin{matrix} \includegraphics[scale=0.6]{trees3.29} \end{matrix}$ is\textit{ not} a Catalan tableaux since it violates the second defining property. 

Catalan tableaux were considered as a special case of permutation tableaux in \cite{perm}, where it is shown that the number of Catalan tableaux of index $n$ is the Catalan number $C_n$. There has been a combinatorial proof of this result by constructing a bijective correspondence with binary trees in \cite{viennot}. Our definition can be seen to be equivalent to the Definition 2.1 of \cite{viennot} by removing the column without $1$'s.

Here we propose another simple bijection, from $\tr_n$ onto the set of Catalan tableaux of index $n + 1$. We shall describe it as a procedure of obtaining a $0/1$ tableaux from a plane binary tree in $3$ steps. As a running example we take $T = 
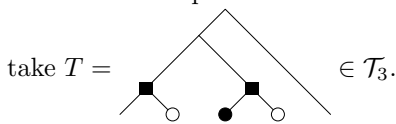
\begin{matrix}
\includegraphics[scale=0.35]{trees3.5}
\end{matrix}   \in \tr_3$. 
\begin{enumerate}
\item We first draw a lattice path in the $2D$ square lattice, starting from anywhere, as we read \textit{all} endpoints of $T$ from left to right: We move one unit to the right(respectively, one unit up) each time we read a left--descendent(right--descendent, respectively). Note that the lattice path has length $n+2$, the first(last) step is always to the right(up, respectively), so that we can form a Young diagram $\lambda(T)$ by adding a vertical segment from the starting point and a horizontal line from the ending point. For example, 
$$\lambda
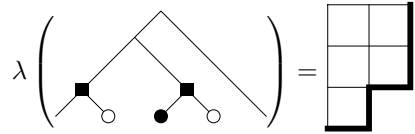
\begin{pmatrix}
\includegraphics[scale=0.35]{trees3.5}
\end{pmatrix} = \begin{matrix}
\includegraphics[scale=0.55]{trees3.25} \end{matrix}.$$
Here the lattice path is drawn in bold. (Remark that $\lambda(T)$ depends only on the reduced configuration $R(T)$)
\item Then we rotate the binary tree by $\pi / 4$ and embed it conformally into the Young diagram  $\lambda(T)$ as illustrated
$$
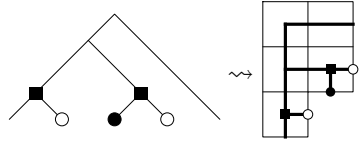
\begin{matrix}
\includegraphics[scale=0.35]{trees3.5}
\end{matrix} \leadsto   \begin{matrix}
\includegraphics[scale=0.6]{trees3.26} \end{matrix} $$
such that each endpoint is at the centre of its corresponding lattice path segment.
\item Finally, to obtain the $0/1$ tableaux, we fill a $1$ into every square inside which the pattern $\begin{matrix}\includegraphics[scale=0.6]{trees3.27}\end{matrix}$ occurs and filling a $0$ into any other square, So the above example gives the tableaux 
$$ 
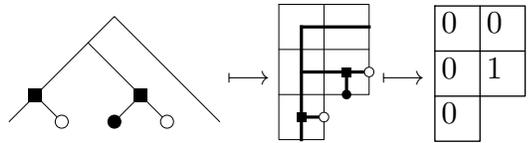
\begin{matrix}
\includegraphics[scale=0.35]{trees3.5}
\end{matrix} \longmapsto   \begin{matrix}
\includegraphics[scale=0.6]{trees3.26} \end{matrix} \longmapsto  \begin{matrix}
\includegraphics[scale=1.2]{trees3.28} \end{matrix}.$$
\end{enumerate}
So far we have defined a map, which we denote $\phi$, from $\tr_n$ to a certain set of $0$/$1$ tableaux. One can show that
\begin{prop}
$\phi$ is a bijection between $\tr_n$ and the set of Catalan tableaux of index $n + 1$. 
\end{prop}

Equipped with this Catalan tableaux $\leftrightarrow$ plane binary tree correspondence, one can reformulate the combinatorial approach to the TASEP stationary state of Corteel and Williams in \cite{corteel4}(more general partially asymmetric processes were considered there using permutation tableaux; one can see that the special case of TASEP corresponds to considering only Catalan tableaux) in terms of planar binary trees: By doing so we recognised that the Markov chain dynamics that the authors of \cite{corteel4} defined on the set of Catalan tableaux can be obtained by using the $\pi$ bijection defined in this note. However, we have refrained from defining a Markov chain dynamics on $\tr_n$ (although it can be done). 

\acknowledgements
The author is indebted to his advisor Kirone Mallick for his patient help during the preparation of this work. He thanks Jeremie Bouttier for valuable suggestions and for fruitful discussions, which have led to the Appendix \ref{app:tableaux} in particular.


\begin{thebibliography}{1}
\bibitem{zeilberger}  L. Shapiro, D. Zeilberger, 1982, \textit{A Markov chain occurring in enzyme kinetics}, J. Math. Biology 15 , 351 –357.
\bibitem{derrida} B. Derrida, E. Domany and D. Mukamel, 1992, \textit{An exact solution of a one-dimensional asymmetric exclusion model with open boundaries,} J. Stat. Phys. 69, 667.
\bibitem{derrida2} B. Derrida, M. R. Evans, V. Hakim, V. Pasquier, 1993,\textit{ Exact solution of a 1D asymmetric exclusion model using a matrix formulation}, J. Phys. A: Math. Gen. 26, 1493.
\bibitem{matrix}  R. A. Blythe and M. R. Evans, 2007, \textit{Nonequilibrium steady states of matrix-product form: a solver’s guide}, J. Phys. A:Math. Theor. 40, R333.
\bibitem{angel} O. Angel,\textit{ The stationary measure of a 2-type totally asymmetric exclusion process}, J. Comb. Theory A 113 , 625 (2006).
\bibitem{duchi} E. Duchi, G. Schaeffer, \textit{A combinatorial approach to jumping particles}, J. Combin. Theory Ser. A 110 (2005), 1–29.
\bibitem{corteel} S.Corteel, R.Brak, A.Rechnitzer,J. Essam, 2005, \textit{A combinatorial derivation of the PASEP stationary state}, FPSAC’05, Taormina. 
\bibitem{corteel2} S. Corteel, L. K. Williams, 2007, \textit{Tableaux combinatorics for the asymmetric exclusion process.} Advances in applied mathematics, 39(3), 293-310.
\bibitem{corteel3} S. Corteel, L. K. Williams, 2011, \textit{Tableaux combinatorics for the asymmetric exclusion process and Askey-Wilson polynomials}, Duke Mathematical Journal, 159(3), 385-415.
\bibitem{corteel4} S. Corteel, L. K. Williams, 2007, \textit{A Markov Chain on Permutations which Projects to the PASEP}, International mathematics research notices, Vol 2007.
\bibitem{mandel}  O. Mandelstam, 2013. \textit{ A Determinantal Formula for Catalan Tableaux and TASEP Probabilities}. arXiv preprint arXiv:1310.7059.
\bibitem{viennot} X.G. Viennot, 2007, \textit{Canopy of binary trees, Catalan tableaux and the asymmetric exclusion process}, in Proc. FPSAC07 (Formal Power Series and Algebraic Combinatorics), Tienjin, Chine, 2007, 12 pp.
\bibitem{2species}  A. Ayyer, J. L. Lebowitz, and E. R. Speer, \textit{On the Two Species Asymmetric Exclusion Process with Semi-Permeable Boundaries}, J. Stat. Phys. 135 
(2009), 1009–1037.
\bibitem{multi} P. A. Ferrari and J. B. Martin,\textit{ Stationary distributions of multi-type totally asymmetric exclusion processes}, Ann. Probab. 35, 807 (2007).
\bibitem{stanley} R. Stanley, \textit{Enumerative Combinatorics,} Vol. 2, Cambridge, 1999.
\bibitem{perm} E. Steingrimsson and L.K. Williams, \textit{Permutation tableaux and permutation patterns}, arXiv:math.CO/0507149, J. Combinatorial Th. A.,114 (2007) 211-234.
\end{thebibliography}
\end{document}